\documentclass[10pt,reqno]{amsart}   
  \usepackage{geometry}
\usepackage{fullpage,amsmath,amsfonts}
\usepackage[all]{xy}
\usepackage{amssymb,amsthm}
\usepackage{hyperref}

\newcommand{\Set}{\mathbf{Set}}
\newcommand{\Sets}{\ensuremath{\Set}}

\newcommand{\calA}{\ensuremath{\mathcal{A}}}
\newcommand{\calB}{\ensuremath{\mathcal{B}}}
\newcommand{\calC}{\ensuremath{\mathcal{C}}}
\newcommand{\calD}{\ensuremath{\mathcal{D}}}
\newcommand{\calE}{\ensuremath{\mathcal{E}}}
\newcommand{\calS}{\ensuremath{\mathcal{S}}}
\newcommand{\calX}{\ensuremath{\mathcal{X}}}

\newcommand{\frakG}{\mathfrak{G}} 

\newcommand{\twopl}[2]{\langle #1, #2\rangle}

\newcommand{\Psh}[1]{\widehat{#1}}
\newcommand{\Shv}{\ensuremath{\mathrm{Sh}}}

\newcommand{\fSet}{\mathbf{fSet}}
\newcommand{\fPos}{\mathbf{fPos}}

\newcommand{\Dec}{\mathrm{Dec}}

\theoremstyle{plain}
\newtheorem{theorem}{Theorem}[section]
\newtheorem{corollary}[theorem]{Corollary}
\newtheorem{proposition}[theorem]{Proposition}
\newtheorem{lemma}[theorem]{Lemma}

\theoremstyle{definition}
\newtheorem{definition}[theorem]{Definition}
\newtheorem{example}[theorem]{Example}
\newtheorem{remark}[theorem]{Remark}

\begin{document}

\title{Decidable objects and Molecular toposes\footnote{Accepted for publication in the {\em Revista de la Uni\'{o}n Matem\'{a}tica Argentina}, \url{https://doi.org/10.33044/revuma}.}}

\author{Mat\'\i as Menni}
\address{Conicet and Universidad Nacional de La Plata, Argentina}
\email{matias.menni@gmail.com}
\urladdr{sites.google.com/site/matiasmenni} 

\begin{abstract}
We study several sufficient conditions for the molecularity/local-connectedness of  geometric morphisms. 
In particuar, we show that if $\calS$ is a Boolean topos then, for every hyperconnected essential geometric morphism ${p : \calE \rightarrow \calS}$ such that the leftmost adjoint $p_!$ preserves finite products, $p$ is molecular and ${p^* : \calS \rightarrow \calE}$ coincides with the full subcategory  of decidable objects in $\calE$. We also characterize the reflections between categories with finite limits that induce molecular maps between the respective presheaf toposes. As a corollary we establish the molecularity of certain geometric morphisms between Gaeta toposes.
\end{abstract}

\maketitle

\tableofcontents

\section{Introduction}

Molecular geometric morphisms were introduced in \cite{BarrPare80} as way to define what it means for an elementary topos to be `locally connected' over another topos.
For instance, for a topological space $X$, the topos of sheaves on $X$ is molecular if and only if $X$ is locally connected. On the other hand, non-localic examples of molecular toposes are the typical   Grothendieck toposes `of spaces' that arise  in Algebraic Geometry \cite{SGA4},  Combinatorial Topology \cite{GabrielZisman}, Synthetic Differential Geometry \cite{KockSDG2ed}, and other areas including  the more recent Rig Geometry sketched in \cite{Lawvere08}.

It is well-known that every essential geometric morphism over $\Set$ is molecular but, apart from the different characterizations of molecular maps \cite[C3.3.1, C3.3.5]{elephant}, there seem to be no practical elementary sufficient conditions. Even over restricted bases. For example, the classical `gros' Zariski topos for an algebraically closed field is well-known to be molecular (over $\Set$) because every essential geometric morphism over $\Set$ is molecular.
On the other hand, for a non-algebraically closed field, the natural base of the corresponding Zariski topos is not $\Set$ but a (Boolean) Galois topos and so, the argument for molecularity must change. Between presheaf toposes, \cite[C3.3.8]{elephant} seems to be the main tool to establish molecularity (and it may be combined with pullback stability to derive molecularity over Grothendieck toposes) but the required hypotheses are far from elementary.

One of the main purposes of this paper is to prove a  sufficient condition for molecularity (Theorem~\ref{Thm}) and discuss some of its implications.
The sufficient  condition states that: if $\calS$ is a Boolean topos then, for every connected essential geometric morphism ${p : \calE \rightarrow \calS}$ such that the leftmost adjoint $p_!$ preserves finite products, $p$ is molecular. To sketch the implications it is convenient to briefly discuss the motivation.

An axiomatic account of   toposes `of spaces'  is proposed in \cite{Lawvere07} via the notion of {\em cohesive} geometric morphism. It is relevant to stress that molecularity is not among the axioms. 
The question if every {\em pre-cohesive} geometric morphism is molecular was left open in \cite{LawvereMenni2015}.
This question stimulated the  recent   \cite{HemelaerRogers2021} and \cite{GarnerStreicher2021} which are devoted to the construction of geometric morphisms that have some of the properties of pre-cohesive maps   but are not molecular. (See also \cite[Remark~{3.8}]{Menni2022}.)  
These constructions might suggest that the question has a negative answer. On the other hand, our sufficient condition implies  that  every pre-cohesive map with Boolean codomain is molecular (Corollary~\ref{CorBetterThanJSL2.6}), which gives a partial positive answer to the question in  \cite{LawvereMenni2015}.  
Moreover, it follows that if $\calS$ is Boolean and ${p : \calE \rightarrow \calS}$ is pre-cohesive then the full inclusion ${p^* : \calS \rightarrow \calE}$ coincides with the subcategory ${\Dec\calE \rightarrow \calE}$ of decidable objects in $\calE$. This improves one of the main results in \cite{Menni2018} and allows us to characterize the quality types  \cite{Lawvere07, MarmolejoMenni2021} over Boolean bases as those pre-cohesive maps whose domain topos is De Morgan (Proposition~\ref{PropQT}).
It also allows us to strengthen one of the main  results  in \cite{LawvereMenni2015}.

The question of what functors between small categories induce molecular maps between the respective toposes does not have a definitive answer but a second purpose  of the paper is to  characterize, in Section~\ref{SecSemiLeftExactness}, the reflections (between categories with finite limits) that do.
This result is of independent interest and it is applied in Section~\ref{SecMolecularBetweenGaeta} to produce molecular maps between Gaeta toposes. Indeed, decidable objects also play a role at the level of sites and, in Section~\ref{SecMolecularBetweenGaeta}, we show that for certain small extensive categories $\calE$, the inclusion ${\Dec\calE \rightarrow \calE}$ induces a molecular geometric morphism between the associated Gaeta toposes. As a byproduct we give a more conceptual account of some constructions in \cite{MenniExercise}.

\begin{remark} 
The first results in the present paper were obtained in late 2020 and the first drafts  were distributed in 2021.
One such draft partially motivated J.~Hemelaer who among other results in \cite{Hemelaer2022} proves (with  very different  techniques) an improvement of  Theorem~\ref{Thm} here.  Some of the ideas there are discussed in Section~\ref{SecBasicToposes}.
\end{remark}

A pleasing aspect of the proofs of some of the main results here is that that they rest on very basic facts about decidable objects in  extensive categories   which seem to have been overlooked until now and that we present in Section~\ref{SecDecidableObjectsAndStableUnits}.

We assume that the reader is familiar with the basic theory of  toposes and geometric morphisms between them (as presented in \cite[Sections~{A2, A4 and C3}]{elephant}, for example) but we recall a fundamental result about extensive categories. 
A  category $\calE$ with finite coproducts is {\em extensive} if the canonical ${\calE/X \times \calE/Y \rightarrow \calE(X+ Y)}$ is an equivalence for every pair of objects $X$, $Y$. See \cite{Gates98a,Carboni96} and references therein.
A map ${X \rightarrow Z}$ in an extensive category is a {\em summand} if there is a map ${Y \rightarrow Z}$ such that the cospan ${X \rightarrow Z \leftarrow Y}$ is a coproduct. It follows that summands are regular monomorphisms. 
Although the definition of extensivity does not require any limits, it follows that summands may be pulled back along any map and that the resulting pulled-back monomorphism is a summand.

\begin{lemma}[Lemma~3.7 in \cite{Gates98a}]
\label{LemGates} For extensive categories $\calE$ and $\calS$ the following hold:
\begin{enumerate}
\item If ${\Psi : \calE \rightarrow \calS}$ is a finite-coproduct preserving functor then $\Psi$ preserves pullbacks of summands.
\item For finite-coproduct preserving functors ${\Psi, \Xi : \calE \rightarrow \calS}$, and any natural transformation ${\Psi \rightarrow \Xi}$, the naturality squares at summands are pullbacks. 
\end{enumerate}
\end{lemma}

\section{Decidable objects and stable units}
\label{SecDecidableObjectsAndStableUnits}

An object $X$  in an extensive category with finite products is  {\em decidable} if the diagonal ${\Delta : X \rightarrow X\times X}$ is a summand \cite{Carboni96,Gates98a}. 
Initial and terminal objects are decidable. 
Despite its simplicity the following result seems to be new.

\begin{proposition}\label{PropMain} Let $\calE$ and $\calS$ be extensive categories with finite products, and let ${\Psi : \calE \rightarrow \calS}$ be a finite-coproduct preserving functor. Then $\Psi$ preserves finite products if and only if it preserves pullbacks over decidable objects.
\end{proposition}
\begin{proof}
One direction is trivial because the terminal object is decidable.
For the other assume that the  square on the left below 
\[\xymatrix{
P \ar[d]_-{\pi_0} \ar[r]^-{\pi_1} & Y \ar[d]^-h &&  P \ar[d]_-{\twopl{\pi_0}{\pi_1}} \ar[r] & S \ar[d]^-{\Delta}  \\
X \ar[r]_-{g}                                 & S                 && X \times Y \ar[r]_-{g\times h} & S \times S
}\]
is a pullback in $\calE$ with decidable $S$. It follows that the square on the right above is also a pullback.
The diagonal of $S$ is complemented so, as $\Psi$ preserves finite coproducts it preserves the pullback on the right above by Lemma~\ref{LemGates}.
Hence, if $\Psi$ also preserves finite products then the square on the left below
\[\xymatrix{
\Psi P \ar[d]_-{\twopl{\Psi \pi_0}{\Psi \pi_1}} \ar[rr] && \Psi S \ar[d]^-{\Delta} && \Psi P \ar[d]_-{\Psi \pi_0} \ar[r]^-{\Psi \pi_1} & \Psi Y \ar[d]^-{\Psi h}  \\
 \Psi X \times \Psi Y \ar[rr]_-{\Psi g\times \Psi h} && \Psi S \times \Psi S && \Psi X \ar[r]_-{\Psi g}                                 & \Psi S 
}\]
is a pullback in $\calS$. In other words, the square on the right above is a pullback.
\end{proof}

If $\calE$ is an extensive category with finite products then we let  ${\Dec \calE \rightarrow \calE}$ be the full subcategory determined by the decidable objects.
This subcategory contains $0$, $1$ and is closed under finite products, finite coproducts and subobjects. 
It follows that $\Dec\calE$ is extensive (and has finite products) \cite{Carboni96}.
The reader is invited to picture the objects of $\calE$ as `spaces' and those  in ${\Dec\calE \rightarrow \calE}$ as `discrete spaces'.
If ${\Dec\calE \rightarrow \calE}$ has a left adjoint ${L : \calE \rightarrow \Dec\calE}$ (and especially if $L$ preserves finite products)  then we picture $L$ as sending a space $X$ to the discrete space ${L X}$ of `connected components' of $X$.

\begin{corollary}\label{CorStableUnits} 
Let $\calE$ be an extensive category with finite limits and assume that the full subcategory ${\Dec\calE \rightarrow \calE}$ has a left adjoint ${L  : \calE \rightarrow \Dec\calE}$. Then $L$ preserves finite products if and only if $L$ preserves pullbacks over decidable objects. 
\end{corollary}
\begin{proof}
Immediate from Proposition~\ref{PropMain} applied to ${L  : \calE \rightarrow \Dec\calE}$.
\end{proof}

The condition that appears in Corollary~\ref{CorStableUnits} has a wider significance that we briefly suggest below in terms of a suggestive reformulation of that result. First  we recall the following concept   \cite[3.7]{CJKP}.

\begin{definition}\label{DefStableUnits}
The adjunction ${L\dashv R}$ with fully faithful $R$ is said to have {\em stable units} if $L$ preserves pullbacks over objects of the form ${R B}$ with $B$ in $\calB$. 
\end{definition}

 Corollary~\ref{CorStableUnits} may then be reformulated as follows:  assuming that ${\Dec\calE \rightarrow \calE}$ is reflective, the left adjoint preserves finite products  if and only if   the adjunction has  stable units.

\begin{example}[Finite posets]\label{ExPosets}
Let ${\fPos}$ be the category of finite posets. It is  not difficult to show that ${\Dec(\fPos) \rightarrow \fPos}$ coincides with the inclusion ${\fSet \rightarrow \fPos}$ of finite discrete posets, and  that this inclusion has a left adjoint (sending a poset to its associated set of  `connected components') that preserves finite products.
\end{example}

\begin{example}[Affine schemes]\label{ExSepAlgs}
Let ${\calE}$ be the opposite of the category of finitely presentable $k$-algebras for a field $k$.
In this case, ${\Dec\calE}$ is the opposite of the category of separable $k$-algebras.
It follows from the results in \cite[I, \S 4, n.6]{DemazureGabriel} that the inclusion ${\Dec\calE\rightarrow \calE}$ has a left adjoint that preserves finite products. 
\end{example}

In an extensive category an object is called {\em connected} if it has exactly two complemented subobjects.
In the two examples above,  the categories involved are essentially small and satisfy that every object is a {\em finite} coproduct of connected objects.
The following two examples are different.

\begin{example}[De Morgan toposes]\label{ExDeMorgan} If $\calE$ is a De Morgan topos then the subcategory ${\Dec \calE \rightarrow \calE}$ coincides with the subcategory of $\neg\neg$-separated objects  by \cite[Proposition~{D4.6.2(iv)}]{elephant}. Subcategories of separated objects for a Lawvere-Tierney topology are well-known to be reflective and such that the left adjoint preserves monomorphisms and finite products  \cite[Theorem~{1.5.1}]{CarboniMantovani}.
\end{example}

\begin{example}[Non-example: topological spaces]\label{ExTop}
In the extensive category of topological spaces and continuous functions between them, an object is decidable if and only if it is discrete.
In contrast with the previous examples, the subcategory of decidable objects  is not reflective. Indeed,  the inclusion does not preserve (infinite) products.
\end{example}

The next  result seems to be new also.
In order to state it  we introduce an ad-hoc piece of terminology.
If ${L\dashv R: \calS \rightarrow \calE}$ is an adjunction  then we say that the  right adjoint $R$ is {\em closed under subobjects} if for every $A$ in $\calS$ and monic ${X \rightarrow R A}$, the unit ${X \rightarrow R(L X)}$ is an isomorphism.

\begin{proposition}\label{PropFPandClosureUnderSubs} Let $\calE$, $\calS$ be extensive categories with finite products and let  ${R : \calS \rightarrow \calE}$  preserve finite coproducts. If $R$ has a finite-product preserving left adjoint  and is closed under subobjects then, for every decidable object $X$ in $\calE$, the unit ${X \rightarrow R(L X)}$ is an isomorphism.
\end{proposition}
\begin{proof}
Let ${L : \calE \rightarrow \calS}$ be the finite-product preserving left adjoint of $R$ and let $\sigma$ be the unit of ${L \dashv R}$.
Notice that both  $L$ and $R$  preserve finite coproducts.
By Lemma~\ref{LemGates}, the unit $\sigma$   is such that the naturality squares at summands are pullbacks.
Let $X$ be a decidable object in $\calE$.
Then the  left inner square below is a pullback
\[\xymatrix{
X \ar[d]_-{\Delta} \ar[r]^-{\sigma} & R(L X) \ar[d]^-{R(L \Delta)} \ar[r]^-{id} & R(L X)  \ar[d]^-{\Delta} \\
X \times X \ar[r]_-{\sigma} \ar@(rd,ld)[rr]_-{\sigma\times\sigma} & R(L  (X \times X))  \ar[r]_-{\cong} & R(L X) \times   R(L X) 
}\]
because the diagonal of $X$ is a summand. As $L$ and $R$ also preserve finite products the right inner square above is a pullback.
Then the rectangle is a pullback, which means that the unit ${\sigma : X \rightarrow R (L X)}$ is monic.
As ${R : \calS \rightarrow \calE}$ is closed under subobjects, ${\sigma_X : X \rightarrow R(L X)}$ is an isomorphism.
\end{proof}

We are mainly interested in cases where $R$ is fully faithful and, for the sake of brevity, we introduce another piece of ad-hoc terminology.

\begin{definition}\label{DefDetached}
 For an extensive category $\calE$ with finite products a  full reflective subcategory ${\calS \rightarrow \calE}$ is {\em detached} if it is  closed under finite coproducts, closed under subobjects, and the left adjoint preserves finite products.
\end{definition}

 For a De Morgan topos $\calE$ as in Example~\ref{ExDeMorgan}, ${\Dec\calE \rightarrow \calE}$ is detached. On the other hand, for an arbitrary extensive $\calE$  with finite products, ${\Dec\calE \rightarrow \calE}$ is closed under finite coproducts, subobjects and finite products but,  it may fail to be reflective, as in Example~\ref{ExTop}. Anyway, the ad-hoc terminology allows us to state the following corollary of Proposition~\ref{PropFPandClosureUnderSubs} succinctly.

\begin{corollary}\label{CorMinimal} Let $\calE$ is an extensive category with finite limits. If  the full subcategory ${\Dec\calE \rightarrow \calE}$ is detached then  it is the least detached subcategory of $\calE$.
\end{corollary}

\section{A sufficient condition for molecularity}
\label{SecSufForMolecularity}

Let ${p : \calE \rightarrow \calS}$ be a geometric morphism. Recall that $p$ is  {\em connected}, if the inverse image ${p^* : \calS \rightarrow \calE}$ of $p$ is fully faithful. Also, $p$ is  {\em essential} if the inverse image ${p^* : \calS \rightarrow \calE}$ of $p$ has left adjoint (usually denoted by ${p_! : \calE \rightarrow \calS}$).

Any functor ${L : \calA \rightarrow \calB}$ determines, for each  object $A$ in $\calA$, the obvious  functor ${L/A :\calA/A \rightarrow \calB/L A}$ that sends ${f : X \rightarrow A}$ to ${L f  :  L X \rightarrow L A}$. If $\calA$ and $\calB$ have finite limits and $L$ has a right adjoint then $L/A$ also has a right adjoint.
In particular, for every geometric morphism ${p : \calE \rightarrow \calS}$ and $A$ in $\calS$ we have a `sliced' geometric morphism ${p/A : \calE/p^* A \rightarrow \calS/A}$ whose inverse image ${(p/A)^* : \calS/A \rightarrow \calE/p^*A}$ is just ${p^*/A}$.

\begin{definition}
A geometric morphism ${p : \calE \rightarrow \calS}$ is {\em molecular} (or {\em locally connected}) if, for each $A$ in $\calS$, the inverse image of   ${p/A : \calE/p^* A \rightarrow \calS/A}$ is cartesian closed.
\end{definition}

Every molecular geometric morphism is essential but the converse does not hold.

\begin{example}\label{ExOfNonMolecularMap}
The essential geometric morphism ${\Set\times \Set \rightarrow \Set^{\rightarrow}}$ whose inverse image sends ${A \rightarrow B}$ to ${(A, B)}$ is not molecular \cite{BarrPare80}.
\end{example}

Recall that for any connected essential geometric morphism ${p : \calE \rightarrow \calS}$,  the leftmost adjoint ${p_! : \calE \rightarrow \calS}$ preserves finite products if and only if  the full subcategory ${p^* : \calS\rightarrow \calE}$ is an exponential ideal   \cite[Proposition~{A4.3.1}]{elephant}.

\begin{theorem}\label{Thm} 
If $\calS$ is a Boolean topos then every connected essential geometric morphism ${p : \calE \rightarrow \calS}$ such that ${p_! : \calE \rightarrow \calS}$ preserves finite products is molecular.
\end{theorem}
\begin{proof}
Proposition~\ref{PropMain} implies that ${p_! : \calE \rightarrow \calS}$ preserves pullbacks over decidable objects.
As $\calS$ is Boolean, every object $A$ in $\calS$ is decidable, so ${p^* A}$ is decidable too.
Then ${p_! : \calE \rightarrow \calS}$ preserves pullbacks over objects of the form ${p^* A}$.
Hence, $p$ is molecular by \cite[Proposition~{10.3}]{LawvereMenni2015}.
\end{proof}

A geometric morphism ${p : \calE \rightarrow \calS}$ is {\em hyperconnected}  if it is connected and the counit of ${p^* \dashv p_*}$ is monic.
For such a geometric morphism the fully faithful   ${p^*}$ is closed under coproducts,  subobjects, and   finite limits, but it may fail to be detached (Definition~\ref{DefDetached}) because it need not have a left adjoint. In any case, with some care, we may also use the intuition that $\calE$ is a category `of spaces' and that  ${p^* : \calS \rightarrow \calE}$ is the full subcategory of `discrete spaces'.
Notice that this is analogous to the intuition that we described for ${\Dec\calE \rightarrow \calE}$.
 The next result shows that the two intuitions are compatible.

\begin{corollary}\label{CorBetterThanJSL2.3} If ${p : \calE \rightarrow \calS}$ is a hyperconnected and essential geometric morphism such that ${p_! : \calE \rightarrow \calS}$ preserves finite products then the  subcategory ${\Dec\calE \rightarrow \calE}$ factors through ${p^* : \calS \rightarrow \calE}$.
\end{corollary}
\begin{proof}
The hypotheses on $p$  imply that ${p^* : \calS\rightarrow \calE}$ is detached so Proposition~\ref{PropFPandClosureUnderSubs}  may be applied to the adjunction ${p_! \dashv p^*}$.
\end{proof}

Corollary~\ref{CorBetterThanJSL2.3} is a strengthening of \cite[Proposition~{2.3}]{Menni2018} and naturally leads to the following strengthening of  \cite[Corollary~{2.4}]{Menni2018}.

\begin{proposition}\label{PropBetterThatJSL2.4} If  $\calS$ is a Boolean topos then, for every hyperconnected essential  geometric morphism ${p : \calE \rightarrow \calS}$  such that ${p_! : \calE \rightarrow \calS}$ preserves finite products,  $p$ is molecular and  ${p^* : \calS \rightarrow \calE}$ coincides with ${\Dec\calE \rightarrow \calE}$.
\end{proposition}
\begin{proof}
Molecularity follows from Theorem~\ref{Thm}.
Also, Corollary~\ref{CorBetterThanJSL2.3} implies that ${\Dec\calE \rightarrow \calE}$ factors through ${p^* : \calS \rightarrow \calE}$.
Finally, since $\calS$ is Boolean, ${p^* : \calS \rightarrow \calE}$ factors through ${\Dec\calE \rightarrow \calE}$ (as in the proof of Theorem~\ref{Thm}).
\end{proof}

\section{Pre-cohesive toposes over Boolean bases}
\label{SecPreCohesiveOverBoolean}

A geometric morphism ${p : \calE \rightarrow \calS}$ is  {\em local} if its direct image  has a fully faithful right adjoint (usually denoted by ${p^! : \calS \rightarrow \calE}$).
Local,   (hyper)connected  and essential geometric morphisms are well-known and  their theory is developed exhaustively  in  \cite[C3]{elephant}.
On the other hand, the following is more recent.

\begin{definition} A  geometric morphism ${p : \calE \rightarrow \calS}$ is {\em pre-cohesive} if it is local, hyperconnected, essential and the leftmost adjoint ${p_!}$  preserves finite products.
\end{definition}

More explicitly, ${p : \calE \rightarrow \calS}$ is pre-cohesive if and only if the adjunction ${p^* \dashv p_*}$ extends to a string ${p_! \dashv p^* \dashv p_* \dashv p^!}$ such that ${p^*, p^! : \calS \rightarrow \calE}$ are fully faithful, the counit of ${p^* \dashv p_*}$ is monic and ${p_! : \calE \rightarrow \calS}$ preserves finite products. The intuition suggested in \cite{Lawvere07} and related references is that  $\calE$ is a category `of spaces', that $\calS$ is a category `of sets'  and that ${p^* : \calS \rightarrow \calE}$ is the full subcategory of `discrete' spaces. So the leftmost adjoint $p_!$ is thought of as a `$\pi_0$' functor assigning to each space the associated set of connected components.
On the other hand, the right adjoint $p_*$ to $p^*$ sends a space to the associated set of points.
Finally, the rightmost adjoint may be pictured as the full subcategory of codiscrete spaces. See \cite{Lawvere07, LawvereMenni2015} and references therein.

Pre-cohesive geometric morphisms were introduced in \cite{MenniExercise} as a weakening of the notion of {\em cohesive} geometric morphism stemming from \cite{Lawvere07}. We stress again, a pre-cohesive geometric morphism is not required to be molecular by definition.
(It is relevant to mention that  a hyperconnected ${p : \calE\rightarrow \calS}$ is pre-cohesive if and only if ${p^* : \calS \rightarrow \calE}$ is cartesian closed and ${p_* : \calE \rightarrow \calS}$ preserves coequalizers \cite[Corollary~{6.2}]{Menni2021}.) 

The following is, at the same time, a partial positive answer to the question about molecularity of pre-cohesive maps, and  a strengthening of \cite[Corollary~{2.6}]{Menni2018}.

\begin{corollary}\label{CorBetterThanJSL2.6} If $\calS$ is Boolean then  every pre-cohesive   ${p : \calE \rightarrow \calS}$ is molecular.
Moreover, ${p^* : \calS \rightarrow \calE}$ coincides with ${\Dec\calE \rightarrow \calE}$ and 
${p^! : \calS \rightarrow \calE}$ coincides with  ${\calE_{\neg\neg} \rightarrow \calE}$.
\end{corollary}
\begin{proof}
Most of the statement follows from  Proposition~\ref{PropBetterThatJSL2.4}. Also,  by \cite[Proposition~{4.4}]{LawvereMenni2015}, the subtopos ${p_* \dashv p^!}$ coincides with ${\calE_{\neg\neg} \rightarrow \calE}$.
\end{proof}

In other words, relying on the terminology used in \cite{Menni2018}:   if $\calS$ is Boolean and    ${p : \calE \rightarrow \calS}$  is pre-cohesive then $p$ is molecular and   ${p_* : \calE \rightarrow \calS}$ is a Unity and Identity for the subcategories ${\Dec\calE \rightarrow \calE}$ and ${\calE_{\neg\neg} \rightarrow \calE}$, making them adjointly opposite.

The axioms for Cohesion are positive and so it is not surprising that any equivalence between toposes is a pre-cohesive geometric morphism.
Still extreme, but more interesting in general, are the pre-cohesive geometric morphisms $p$ such that ${p_! = p_*}$; equivalently, such that ${p_* \dashv p^*}$. These are called {\em quality types}, see \cite{Lawvere07} and also \cite{MarmolejoMenni2021}.
Over a Boolean base $\calS$, Corollary~\ref{CorBetterThanJSL2.6} leads to the following characterization of quality types ${p : \calE \rightarrow \calS}$ in terms of the internal logic of $\calE$.

\begin{proposition}\label{PropQT} If $\calS$ is a Boolean topos and  ${p : \calE \rightarrow \calS}$ is a  pre-cohesive geometric morphism then, 
 $p$ is a quality type if and only if $\calE$ is De Morgan.
\end{proposition}
\begin{proof}
If $p$ is a quality type then $\calE$ is De Morgan by \cite[Lemma~{5.3(3)}]{MarmolejoMenni2021}.
For the converse first recall that, as observed in \cite{Lawvere07} (see also \cite{Johnstone2011}) there is, for any pre-cohesive geometric morphism $p$,  a canonical monic comparison ${\phi : p^* \rightarrow p^!}$ defined as follows
\[\xymatrix{
p^* \ar[rr]^-{p^* \epsilon^{-1}} && p^* p_* p^! \ar[rr]^-{\beta_{p^!}} && p^!
}\]
where $\epsilon$ is the (necessarily iso)  counit of ${p_* \dashv p^!}$ and $\beta$ is the (monic) counit of ${p^* \dashv p_*}$.
Notice that $p$ is a quality type if and only if $\phi$ is an isomorphism.

To complete the proof assume that ${p : \calE \rightarrow \calS}$ is pre-cohesive with $\calE$ De Morgan and $\calS$ Boolean.
Corollary~\ref{CorBetterThanJSL2.6} implies that  ${p^* : \calS \rightarrow \calE}$ coincides with the inclusion  ${\Dec\calE \rightarrow \calE}$ and ${p^! : \calS \rightarrow \calE}$ coincides with the full subcategory of $\neg\neg$-sheaves.
Since $\calE$ is De~Morgan, ${\Dec\calE \rightarrow \calE}$ coincides with the full subcategory of $\neg\neg$-separated objects (Example~\ref{ExDeMorgan}) and, therefore, ${p^!}$ factors through ${p^*}$. In other words, ${\beta_{p^!} : p^* p_* p^! \rightarrow p^!}$ is an isomorphism and hence $\phi$ is an isomorphism.
\end{proof}

Some of the main results in \cite{LawvereMenni2015} may also be improved in several ways using the sufficient condition for molecularity proved here. 
We do not discuss these in detail but we just sketch one such result to illustrate the idea. We assume familiarity with that paper but we recall some of the terminology and results there.

If ${p : \calE \rightarrow \calS}$ is pre-cohesive then an object $X$ in $\calE$ is {\em connected} if ${p_! X = 1}$.
Intuitively, $X$ is connected if its discrete space of connected components is terminal.
For example, \cite[Proposition~4]{Lawvere07} implies that the subobject classifier $\Omega$ of $\calE$ is connected if and only if all injective objects in $\calE$ are connected. If the base topos $\calS$ is De Morgan then connectedness of $\Omega$ has an alternative formulation that we discuss next.

In any topos we let $2$ denote the coproduct ${1 + 1}$ of the terminal object with itself.
If ${p : \calE \rightarrow \calS}$ is pre-cohesive then intuition about $p^!$ suggests that we may picture ${p^! 2}$ as the codiscrete interval determined by two points. We say that $p$ has the {\em Connected Interval} (CI) property if ${p^! 2}$ is connected.
For instance, the Corollary in \cite[Section~{VI}]{Lawvere07} implies that: if $\calS$ is De Morgan, then $\Omega$ in $\calE$  is connected if and only if $p$ has CI. Hence, combining Proposition~\ref{PropQT} and \cite[Proposition~4]{Lawvere07} we may conclude that: if $\calS$ is Boolean, $\calE$ is De Morgan and $p$ has CI then $\calS$ is inconsistent.

Notice that the CI property makes sense  for any local and essential geometric morphism.
For instance we reiterate that, for pre-cohesive ${p : \calE \rightarrow \calS}$, we still don't know if the sliced ${p/B}$  is pre-cohesive for each $B$ in $\calS$, but we do know that ${p/B}$ is local, hyperconnected and essential (see the discussion before \cite[Definition~{5.7}]{LawvereMenni2015}).

\begin{lemma}\label{LemCIiffStableCI} 
If ${p : \calE \rightarrow \calS}$ is pre-cohesive  then: $p$ has CI if and only if, for each $B$ in $\calS$, ${p/B}$ has CI. 
\end{lemma}
\begin{proof}
This is analogous to \cite[Lemma~{5.8}]{Lawvere07} but apparently stronger in view of \cite[Corollary~{6.5}]{Lawvere07}.
Anyway, one direction is immediate. For the other assume that $p$ has CI.

As we have already recalled, $p/B$ is essential and local. 
Indeed,  for each $B$ in $\calS$ the rightmost adjoint ${(p/B)^! : \calS/B \rightarrow \calE/p^* B}$ may be described as follows \cite[Lemma~{5.4}]{LawvereMenni2015}.
For every ${f : A \rightarrow B}$ in $\calS$, if we let the square on the left below 
\[\xymatrix{
\Phi f \ar[d]_-{f_0} \ar[r]^-{f_1} & p^* B \ar[d]^-{\phi}  && p^! A \times p^* B \ \ar[d]_-{id \times \phi} \ar[r]^-{\pi_1} & p^* B \ar[d]^-{\phi} \\
p^! A \ar[r]_-{p^! f} & p^! B && p^! A \times p^! B \ar[r]_-{\pi_1} & p^! B
}\]
 be a pullback  then ${(p/B)^! f = f_1}$. In particular, for ${f = \pi_1 : A \times B \rightarrow B}$ we obtain the pullback on the right above.
Hence,  
\[ (p/B)^! (\pi_1 : A\times B \rightarrow B) = (\pi_1  : p^! A \times p^* B \rightarrow p^* B)\]
 in $\calE/p^*B$ for any $A$, $B$ in $\calS$.

For $B$ in $\calS$,  the object $2$ in $\calS/B$ is, as a map in $\calS$, the projection ${\pi_1 : 2 \times B \rightarrow B}$,  so  the previous paragraph  implies that ${(p/B)^! 2 = (\pi_1 : p^! 2 \times p^* B \rightarrow B)}$.
By  \cite[Lemma~{5.2}]{LawvereMenni2015}, the leftmost adjoint ${(p/B)_! : \calE/p^* B \rightarrow \calS/B}$ sends ${(p/B)^! 2}$ to the composite on the left below
\[\xymatrix{
p_! ( p^! 2 \times p^* B ) \ar[r]^-{p_! \pi_1}  & p_!(p^* B) \ar[rr]^-{\textnormal{counit}} && B &&  p_!(p^! 2) \times B \ar[r]^-{\pi_1} & B
}\]
but, as $p_!$ preserves finite products, this is isomorphic to the composite on the right above, as  objects in ${\calS/B}$.
In turn, as $p$ has CI, the right map above  is isomorphic (over $B$) to the identity on $B$. 
Altogether, ${(p/B)_! ((p/B)^! 2) = 1}$ in $\calS/B$.
\end{proof}

If ${p : \calE \rightarrow \calS}$ is an essential geometric morphism such that $p_!$ preserves finite products then, for each $X$ in $\calE$ and $S$ in $\calS$, the composite
\[\xymatrix{
p_!(X^{p^* S}) \times S \ar[rr]^-{id \times \tau^{-1}} && p_!(X^{p^* S}) \times p_!(p^* S) \ar[r]^-{\cong} & p_!(X^{p^* S} \times p^* S) \ar[r]^-{p_! ev} & p_! X
}\]
 (where $\tau$ is the counit of ${p_! \dashv p^*}$) transposes to a map ${p_!(X^{p^* S}) \rightarrow (p_! X)^S}$. The map    ${p : \calE \rightarrow \calS}$ is said to satisfy {\em Continuity} if the canonical  ${p_!(X^{p^* S}) \rightarrow (p_! X)^S}$ is an isomorphism   for every $X$ in $\calE$ and $S$ in $\calS$.
A geometric morphism is {\em cohesive} if it is pre-cohesive and satisfies Continuity.

\begin{proposition}\label{PropCohesive}
For any cohesive ${p : \calE \rightarrow \calS}$ satisfying CI,  $p$ is molecular if and only if $\calS$ is Boolean.
\end{proposition}
\begin{proof}
If $\calS$ is Boolean then $p$ is molecular by Corollary~\ref{CorBetterThanJSL2.6}.
Conversely, if $p$ is molecular then it is {\em stably} pre-cohesive by \cite[Corollary~{10.4}]{LawvereMenni2015} in the sense that ${p/B : \calE/p^* B \rightarrow \calS/B}$ is pre-cohesive for every $B$ in $\calS$.
For each such $B$, ${ \calE/p^* B}$ satisfies CI by Lemma~\ref{LemCIiffStableCI} and, since ${(p/B)_!}$ preserves finite products, $p/B$ satisfies Connected Codiscreteness (CC) by \cite[Lemma~{6.2}]{LawvereMenni2015}.
In other words, $p$ satisfies Stable Connected Codiscreteness (SCC); so $\calS$ satisfies the Internal Axiom of Choice by \cite[Corollary~{9.4}]{LawvereMenni2015}, and hence $\calS$ is Boolean \cite[Remark~{D4.5.8}]{elephant}.
\end{proof}

Notice that the proof also shows that if the equivalent conditions of Proposition~\ref{PropCohesive} hold then the internal Axiom of Choice holds in $\calS$.

\section{Molecularity and semi-left-exactness}
\label{SecSemiLeftExactness}

In this section we prove a characterization of reflections (between small categories with finite limits)  that induce molecular maps between the respective presheaf toposes. This result will be applied later to prove molecularity of morphisms between Gaeta toposes.
We assume that the reader is familiar with `Street fibrations'; in particular, with the fact that any such is isomorphic to  a composite of an equivalence and a  `Grothendieck' fibration \cite[(5.1)]{Street1980}. 

Let $\calX$ and $\calB$ be categories with finite limits.
Let ${L : \calX \rightarrow \calB}$ be a functor with fully faithful right adjoint $R$, and let $\sigma$ be the unit of the adjunction.

\begin{definition}
The adjunction ${L\dashv R}$ is said to be {\em semi-left-exact} if for every pullback
\[\xymatrix{
Y  \ar[d] \ar[r]^-u & R B \ar[d]^-g \\
X \ar[r]_-{\sigma} & R(L X)
}\]
with $X$ in $\calX$ and $B$ in $\calB$, ${L u}$ is an isomorphism. 
\end{definition}

It is well-known that stability of units implies semi-left-exactness \cite[(3.7)]{CJKP}.

The following result is folklore and I learned the proof below from G.~Janelidze.

\begin{lemma}\label{LemJanelidze} With  ${L \dashv R}$ as above, $L$ is a Street fibration if and only if the adjunction ${L \dashv R}$ is semi-left-exact.
\end{lemma}
\begin{proof}
For each $A$ in $\calA$, the functor ${L/A : \calA/A \rightarrow \calB/L A}$ that sends ${(f : X  \rightarrow A)}$ to ${L f : L X \rightarrow L A}$ has a right adjoint that we denote by $R_A$.
The straightforward (folk?) variant of the equivalence between the first two items of  \cite[Theorem~{2.10}]{Gray1966} says  that: 
$L$ is a Street fibration if and only if ${R_A :\calB/L A \rightarrow \calA/A}$ is fully faithful for every $A$ in $\calA$.
The second paragraph of \cite[3.6]{CJKP} shows that fully faithfullness of the $R_A$'s is equivalent to semi-left-exactness.
\end{proof}

The adjunction ${L \dashv R : \calB \rightarrow \calX}$ determines an essential and  local geometric morphism ${L : \Psh{\calX} \rightarrow \Psh{\calB}}$ with ${L_! : \Psh{\calX} \rightarrow \Psh{\calB}}$ the left Kan extension  of $L$.

\begin{theorem}\label{ThmMainSMLE}
With  ${L \dashv R : \calB \rightarrow \calX}$ as above, the following are equivalent: 
\begin{enumerate}
\item The geometric morphism ${L : \Psh{\calX} \rightarrow \Psh{\calB}}$  is molecular.
\item The adjunction  ${L \dashv R : \calB \rightarrow \calX}$ is semi-left-exact.
\end{enumerate}

\end{theorem}
\begin{proof}
The first  item implies the second by \cite[Lemma~{3.1(iii)}]{GarnerStreicher2021}.
The converse follows from \cite[Proposition~{C3.3.12}]{elephant}. We give some of the details, trying to be consistent with the notation in that result.
By Lemma~\ref{LemJanelidze} we may assume that $L$ is a fibration and so ${(L, R) : (\calX, M_{\calX}) \rightarrow (\calB, M_{\calB})}$ is a fibration of sites \cite[Definition~{C2.5.6}]{elephant}, where ${M_{\calX}}$ and ${M_{\calB}}$ are the respective minimal Grothendieck coverages (i.e., only maximal sieves cover).
So, in order to apply \cite[Proposition~{C3.3.12}]{elephant}, we need only check that maximal sieves in $\calX$ are $M_{\calB}$-locally connected in the sense of that result. To do this let $V$ be an object of $\calX$, let  ${b_1 : V_1 \rightarrow V}$ and ${b_2 : V_2 \rightarrow V}$ be two maps (i.e. two maps in the maximal sieve on $V$) and let the following square
\[\xymatrix{
U \ar[d]_-{a_1} \ar[r]^-{a_2} & L V_2 \ar[d]^-{L b_2} \\
L V_1 \ar[r]_-{L b_1} & L V
}\]
commute in $\calB$. Trivially, we can connect $b_1$ and $b_2$ by a zigzag of morphisms in  the maximal sieve on $V$ as on the left below
\[\xymatrix{
V_1 \ar[rd]_-{b_1} \ar[r]^-{b_1} & V \ar[d]^-{id} & \ar[l]_-{b_2} \ar[dl]^-{b_2} V_2 && 
  & \ar[ld]_-{a_1} U  \ar[d] \ar[rd]^-{a_2} \\
 & V & && L V_1 \ar[r]_-{L b_1} & L V & L V_2 \ar[l]^-{L b_2}
}\]
and we can find a cone over the image of this zigzag under $L$ with vertex $U$, as on the right above, connecting ${a_1}$ and ${a_2}$.
So, it easily follows that if we take the maximal sieve on $U$ then, for every map in it, that is, for every ${U' \rightarrow U}$ in $\calB$, we can find a cone over the image of the zigzag connecting $b_1$ and $b_2$, connecting the morphisms ${U' \rightarrow U \rightarrow L V_1}$ and ${U' \rightarrow U \rightarrow L V_2}$. 
\end{proof}

Theorem~\ref{ThmMainSMLE} deserves a direct proof not involving Grothendieck topologies.

\begin{corollary}\label{CorMain}
With  ${L \dashv R : \calB \rightarrow \calX}$ as above, the following are equivalent: 
\begin{enumerate}
\item The geometric morphism ${L : \Psh{\calX} \rightarrow \Psh{\calB}}$  is molecular and $L_!$ preserves finite products.
\item The adjunction ${L_! \dashv L^*}$ has stable units.
\item The adjunction  ${L \dashv R}$ is semi-left-exact and $L$ preserves finite products. 
\end{enumerate}
\end{corollary}
\begin{proof}
The first and second items are equivalent by \cite[Proposition~{10.3}]{LawvereMenni2015}.
The first and third items are equivalent by  Theorem~\ref{ThmMainSMLE} together with the well-known fact that  a functor preserves finite products if and only if its left Kan extension does  \cite{BorceuxDay1977}.
\end{proof}

\section{Molecular maps between Gaeta toposes}
\label{SecMolecularBetweenGaeta}

The role of Gaeta toposes as paradigmatic examples of toposes `of spaces' is discussed in   \cite{Lawvere91}; see also \cite[Section~5]{Lawvere08}.
In this section we combine the results in previous sections to give a sufficient condition for an inclusion  ${\Dec\calE \rightarrow \calE}$, for small extensive $\calE$ with finite limits, to induce a molecular (essential and local) map between the associated Gaeta toposes. 
We also explain how to restrict these to (molecular and) pre-cohesive geometric morphisms.
As usual, we denote the topos of presheaves on $\calC$ by $\Psh{\calC}$.

\begin{corollary}\label{CorPreGaeta} If  $\calE$ is a small extensive category with finite limits then restriction along ${\Dec\calE\rightarrow \calE}$ is the direct image of a local geometric morphism ${\Psh{\calE} \rightarrow \Psh{\Dec\calE}}$. If, moreover, the inclusion  ${\Dec\calE\rightarrow \calE}$ has a finite-product preserving left adjoint then the local ${\Psh{\calE} \rightarrow \Psh{\Dec\calE}}$ is also molecular and the leftmost adjoint preserves finite products.
\end{corollary}
\begin{proof}
The first part follows as in the comment before Proposition~\ref{PropMain}.
In more detail, as ${\Dec\calE \rightarrow \calE}$ preserves finite limits, restriction along it determines the direct image of a geometric morphism ${\Psh{\calE} \rightarrow \Psh{\Dec\calE}}$ by \cite[A4.1.10]{elephant}. The direct image has a right adjoint (by the existence of right Kan extensions) which is fully faithful because ${\Dec\calE \rightarrow \calE}$ is  (see, e.g., \cite[VII.4]{maclane2}).
 In other words, the geometric ${f : \Psh{\calE} \rightarrow \Psh{\Dec\calE}}$ is local.

Let ${L : \calE \rightarrow \Dec\calE}$ be the left adjoint to the inclusion in the opposite direction.
The functor $L$  preserves pullbacks over decidable objects by Corollary~\ref{CorStableUnits} so $L$ and its right adjoint form an adjunction with stable units (and hence, semi-left-exact). Corollary~\ref{CorMain} completes the proof.
\end{proof}

Let $\calC$ be an extensive category.
For each object $C$ in $\calC$ we let ${K_{\calC} C = K C}$ be the collection of finite families ${(C_i \rightarrow C \mid i \in I)}$ such that the induced map ${\sum_{i\in I} C_i \rightarrow C}$ is an isomorphism.
It is easy to check that $K$ satisfies axioms for bases of Grothendieck topologies so, if $\calC$ is small, we obtain a Grothendieck topos ${\frakG \calC = \Shv(\calC, K)}$, sometimes called the {\em Gaeta topos} (of $\calC$). See, for example, \cite{Lawvere91} or \cite[Section~2]{Lawvere08}.

Let ${\calD}$ be another extensive category and let ${F : \calC \rightarrow \calD}$ be a functor.
If $F$ preserves finite coproducts then it is clear that, for any ${(C_i \rightarrow C \mid i \in I)}$ in ${ K_{\calC} C}$, the induced family ${(F C_i \rightarrow F C \mid i \in I)}$ is in ${ K_{\calD} (F D)}$.

So, if $\calC$ and $\calD$ are also small and have finite limits and moreover $F$ preserves finite limits then ${F}$ is a morphism of sites ${(\calC, K_{\calC}) \rightarrow (\calD, K_{\calD})}$ in the sense of \cite[Definition~{C2.3.1}]{elephant} and so, by \cite[Corollary~{C2.3.4}]{elephant}, the functor ${F^* : \Psh{\calD} \rightarrow \Psh{\calC}}$ restricts to a functor  ${\frakG \calD \rightarrow \frakG\calC}$  which is the direct image of a geometric morphism.

\begin{lemma}\label{LemGaeta} Let $\calC$ and $\calD$ be small extensive categories with finite limits.
For any functor ${F : \calC \rightarrow \calD}$ preserving finite limits and finite coproducts, the following diagram is a pullback in the category of toposes
\[\xymatrix{
\frakG \calD \ar[d] \ar[r] & \Psh{\calD} \ar[d] \\
\frakG \calC \ar[r] & \Psh{\calC}
}\]
where the horizontal maps are the obvious sutoposes and the vertical maps are the geometric morphisms whose direct image is ${F^*}$.
\end{lemma}
\begin{proof}
Any finite family  ${(D_i \rightarrow D \mid i \in I)}$ in ${K_{\calD} D}$ is the obvious pullback of the family ${(in_i : 1 \rightarrow I\cdot 1 \mid i\in I )}$ in ${K_{\calD} (I\cdot 1)}$ which is, in turn, the image under $F$ of the family  ${(in_i : 1 \rightarrow I\cdot 1 \mid i\in I )}$ in ${K_{\calC} (I\cdot 1)}$.
It follows that ${K_{\calD}}$ is the smallest basis for a Grothendieck topology containing the families ${(F C_i \rightarrow F C \mid i \in I)}$ induced by the  families ${(C_i \rightarrow C \mid i \in I)}$ in ${ K_{\calC} C}$.
The result then follows from \cite[Lemma~{C2.3.12}]{elephant}.
\end{proof}

In particular, for a small extensive category  with finite limits, Lemma~\ref{LemGaeta} applies to the full subcategory of decidable objects.

\begin{lemma}\label{LemLocalBetweenGaeta}
If $\calE$ is a small extensive category with finite products then pre-composition with ${\Dec\calE \rightarrow \calE}$ restricts to the direct image of a local geometric morphism ${\frakG\calE \rightarrow \frakG(\Dec\calE)}$.
\end{lemma}
\begin{proof}
Combine Lemma~\ref{LemGaeta} with the fact that local maps are stable under pullback \cite[C3.6.7(iv)]{elephant} so  $k$ is local.
\end{proof}

\begin{example} If we let $\calE$ be the category of finite posets as in Example~\ref{ExPosets}, then the resulting geometric morphism ${\frakG \calE \rightarrow \frakG(\Dec\calE) \cong \Sets}$ is not only local but it is well-known to be molecular and  pre-cohesive (and that it  classifies distributive lattices with exactly two complemented elements).
\end{example}

\begin{example} If we let $\calE$ be the opposite of the category of finitely presentable $k$-algebras as in  Example~\ref{ExSepAlgs} then ${\frakG \calE \rightarrow \frakG(\Dec\calE) = \calS}$ is local, molecular and the leftmost adjoint preserves finite products by Theorem~\ref{ThmGaeta}, but it is not pre-cohesive if $k$ is not algebraically closed.
\end{example}

\begin{remark}\label{RemNonACfields}
To understand what fails for non algebraically-closed fields it is  convenient to argue in more generality. Let $\calA$  be an essentially small category and let ${L \dashv R : \calB \rightarrow \calA}$ be a reflective subcategory with unit denoted by ${\nu}$. Let ${f : \Psh{\calA} \rightarrow \Psh{\calB}}$ be the induced essential and local geometric morphism and let $\sigma$ be the unit of ${f_! \dashv f^*}$. For representable ${X = \calA(-,A)}$ in $\Psh{\calA}$, 
\[\xymatrix{
 X = \calA(-,A) \ar[rrr]^-{\sigma_X = \calA(-,\nu_A)} &&&   f^*(f_! X) =  \calB(L -, L A) \cong \calA(-, R (L A)) 
}\]
because the left Kan extension $f_!$ of $L$ preserves representables. So, if $\nu_A$ does not have a section then $\sigma_X$ is not epic and, hence, $\sigma$ is not epic. In other words, the Nullstellensatz condition in the definition of pre-cohesive map does not hold; equivalently, $f$ is not hyperconnected. In the case of affine $k$-schemes, this is about the existence of one such scheme $A$ whose decidable reflection ${A \rightarrow \pi_0 A}$ does not have a section; dually, a finitely generated $k$-algebra whose largest separable $k$-subalgebra does not have a retraction.
\end{remark}

\begin{theorem}\label{ThmGaeta} 
Let $\calE$ be a small extensive category with finite limits. If the inclusion  ${\Dec\calE\rightarrow \calE}$ has a finite-product preserving left adjoint then the local map  ${\frakG\calE \rightarrow \frakG(\Dec\calE)}$ (Lemma~\ref{LemLocalBetweenGaeta})
is  molecular and the leftmost adjoint preserves finite products.
\end{theorem}
\begin{proof}
Lemma~\ref{LemGaeta} implies that the following diagram is a pullback
\[\xymatrix{
\frakG\calE \ar[d]_-k \ar[r]^-h & \Psh{\calE} \ar[d]^-f \\
\frakG(\Dec\calE) \ar[r]_-g & \Psh{\Dec\calE}
}\]
where $g$ and $h$ are the evident subtoposes and $f$ is the local geometric morphism induced by ${\Dec\calE \rightarrow \calE}$.
If the inclusion ${\Dec\calE \rightarrow \calE}$ has a finite-product preserving left adjoint then, by Corollary~\ref{CorPreGaeta}, $f$ is molecular and $f_!$ preserves finite products. 
By \cite[Theorem~{C3.3.15}]{elephant}, $k$ is molecular and also the Beck-Chevalley natural transformation ${f^* g_* \rightarrow h_* k^*}$ is an isomorphism. 
It follows that  $k^*$ has a left  adjoint of ${k_! : \frakG\calE \rightarrow \frakG(\Dec\calE)}$ which may be identified with the composite
\[\xymatrix{
\frakG\calE  \ar[r]^-{h_*} &  \Psh{\calE} \ar[r]^-{f_!} & \Psh{\Dec\calE} \ar[r]^-{g^*} & \frakG(\Dec\calE) 
}\]
and, since $f_!$ preserves finite products, so does $k_!$. 
\end{proof}

If  $k$ is a perfect field, not necessarily algebraically closed, \cite[Proposition~{4.3}]{MenniExercise} shows that, pulling the above geometric morphism ${\frakG \calE \rightarrow \frakG(\Dec\calE) = \calS}$  along the subtopos ${\calS_{\neg\neg} \rightarrow \calS}$ results in a pre-cohesive (and molecular) geometric morphism ${\mathcal{F} \rightarrow \calS_{\neg\neg}}$.
Part of the argument may be generalized as follows.

First recall that for a small category $\calB$ equipped with a Grothendieck topology $J$, a map ${\lambda : F \rightarrow G}$ in $\Psh{\calB}$ is {\em locally surjective (w.r.t. $J$)} \cite[Corollary~{III.7.5}]{maclane2} if for each object $B$ of $\calB$ and each 
element ${y\in G B}$, there is a $J$-cover $S$ of $B$  such that for all ${f: C \rightarrow B}$ 
in $S$ the element ${y\cdot f }$  is in the image of ${\lambda_C :  F C \rightarrow G C}$.

Let ${L \dashv R : \calB \rightarrow \calA}$ be an adjunction between small categories and let ${\nu}$ be its unit.
If $\calB$ is equipped with a Grothendieck topology $J$ then, for each $C$ in $\calA$, we may ask if the  natural transformation ${\calA(R-, \nu_C): \calA(R-,C) \rightarrow \calA(R-, R(L C))}$ in $\Psh{\calB}$  is locally surjective; that is, if for every $D$ in $\calB$ and   every ${d : R D \rightarrow R(L C)}$ in $\calA$ there exists a covering sieve ${S \in J D}$ such that, for every ${e : E \rightarrow D}$ in $S$, ${d (R e) : R E \rightarrow R(L C)}$ factors through ${\nu : C \rightarrow R(L C)}$.

\begin{lemma}\label{LemLocalSurjectivity} 
With ${L\dashv R}$, $\nu$ and $J$ as above, let $R$ be fully faithful and   ${f : \Psh{\calA} \rightarrow \Psh{\calB}}$  be the local and essential geometric morphism induced by the adjunction ${L \dashv R}$. The canonical transformation ${f_* X \rightarrow f_! X}$ is locally surjective for every $X$ in $\Psh{\calA}$ if and only if 
${\calA(R-, \nu_C): \calA(R-,C) \rightarrow \calA(R-, R(L C))}$ is locally surjective for every $C$ in $\calA$.
\end{lemma}
\begin{proof}
This is a variant of  \cite[Lemma~{1.15}]{MenniExercise} but we sketch some of the details.
By the results in \cite{Johnstone2011}, the canonical ${\theta_ X : f_* X \rightarrow f_! X}$ is locally surjective if and only if ${f_* \sigma_X : f_* X \rightarrow f_* (f^* (f_! X))}$ is locally surjective where $\sigma$ is the unit of ${f_! \dashv f^*}$. 

If ${f_* \sigma_X}$ is locally surjective for every $X$ then, taking   ${X = \calA(-,C)}$ we may conclude  that  ${f_* \sigma_X = \calA(R-, \nu_C): \calA(R-,C) \rightarrow \calA(R-, R(L C))}$ is locally surjective. (Bear in mind Remark~\ref{RemNonACfields} for the relation between $\sigma$ and $\nu$.)

For the converse let $X$ in $\Psh{\calA}$ and $B$ in $\calB$. Recall that there is a quotient  
\[\xymatrix{ \sum_C X C \times \calB(B, L C) \ar[r] &  (f_! X) B }
 \]
 and that, for ${x \in X C}$ and ${b : B \rightarrow L C}$ in $\calB$, the image of ${(C, x, b)}$ under the quotient above  is denoted by ${x\otimes b \in (f_! X) B  }$.
 The quotient ensures that 
 \[ (x \cdot u)\otimes v = x \otimes ((L u) v) \]
  for every map ${u : D \rightarrow C}$ in $\calA$ and every map  ${v : B \rightarrow L D}$ in $\calB$.
 
 If we let $\xi$ be the (iso) counit of the adjunction ${L \dashv R}$ then the canonical  ${\theta_{X, B} : (f_* X) B \rightarrow (f_! X) B}$  sends ${x \in (f^* X)B = X(R B)}$ to ${x \otimes \xi_B^{-1} \in (f_! X)B}$; notice that ${\xi_B : L (R B) \rightarrow B}$.

 We must show that ${\theta_X :  f_* X \rightarrow f_! X}$ is locally surjective.
 For  ${x\otimes b \in (f_! X)B}$,  consider ${R b \in \calA(R B, R(L C)) = (\calA(R- , R(L C))) B}$.
 By hypothesis,   there exists a $J$-cover ${S \in  J  B}$ such that for every ${e : E \rightarrow  B}$ in $S$, ${(R b) \cdot e = (R b) (R e)}$ is in the image of ${f^* \sigma)}$; in other words,  there exists an  ${e' : R E \rightarrow C}$ in $\calA$ such that ${(f_* \sigma) e' = \nu_C e' = (R b) (R e) = R(b e)}$.
 
We prove that the same $S$ witnesses that ${x\otimes b \in (f_! X)B}$ is locally in the image of ${\theta_X : f_* X \rightarrow f_! X}$.
For each ${e : E \rightarrow  B}$ in $S$, 
\[ b e \xi_E = \xi_{L C} L ((R b) (R e)) = \xi_{L C} (L \nu_C) (L e') = L e' \]
so ${b e = (L e') \xi_E^{-1}}$ and then 
\[ \theta(x\cdot e') = (x\cdot e')\otimes \xi^{-1} = x \otimes ((L e') \xi^{-1}) = x\otimes (b e) = (x\otimes b)\cdot e \]
showing that ${(x\otimes b)\cdot e \in (f_! X)E}$ is in the image of $\theta_X$.
\end{proof}

The intuition is that local surjectivity of ${\calA(R-, \nu_C)}$ for every $C$ is a form of `Nullstellensatz'  for the adjunction ${L\dashv R}$ together with the topology $J$ on $\calB$.

\begin{example} Let $\calA$ have a terminal object and let ${R : \calB \rightarrow \calA}$ be the full (reflective)  degenerate subcategory consisting only of the terminal object. If we endow $\calB$ with the  extreme topology $J$ such that only the maximal sieve covers then, for $C$ in $\calA$, ${\calA(R-, \nu_C): \calA(R-,C) \rightarrow \calA(R-, R(L C))}$ is locally surjective if and only if the object $C$ has a point.
(The  coincidence with the condition appearing \cite{Johnstone2011} is, of course, no accident.)
\end{example}

\begin{example}\label{ExPosetsNew} Let $\calA$ be the category of finite partially ordered sets, equip ${\calB = \Dec\calA = \fSet}$ with the Gaeta topology and consider the reflective subcategory ${\fSet \cong \Dec\calA \rightarrow \calA}$. Then  ${\calA(R-, \nu_C): \calA(R-,C) \rightarrow \calA(R-, R(L C))}$ is locally surjective for every $C$ in $\calA$ because every non-initial poset has a point and the initial set is covered by the empty family.
\end{example}

See also \cite{Menni2021a} for other examples, similar to Example~\ref{ExPosetsNew} and to complex affine schemes, but arising from categories of rigs with idempotent addition.

If every object in $\calA$ is a finite coproduct of connected objects then it is natural to concentrate on the latter. 
Notice that, in this case, the Comparison Lemma implies that $\frakG\calA$ is equivalent to the topos of presheaves on the category of connected objects in $\calA$.
This comment applies to the example above and also to the case where $\calA$ is the category of `affine schemes' for a field $k$, as in the next example.

\begin{example} Let $k$ be a perfect field and let $\calA$ be the opposite of the category of  finitely presentable  $k$-algebras without idempotents.
Let ${R : \calB \rightarrow \calA}$ be the full (reflective) subcategory determined by the separable $k$-algebras without idempotents. If we equip $\calB$ with the atomic topology then Hilbert's Nullstellensatz implies that  the transformation ${\calA(R-, \nu_C): \calA(R-,C) \rightarrow \calA(R-, R(L C))}$ is locally surjective for every $C$ in $\calA$. This is just a reformulation of the main examples in \cite{MenniExercise}.
\end{example}

We next state an analogue of \cite[Proposition~{1.16}]{MenniExercise}, with essentially the same proof.

\begin{proposition}\label{PropStablyPreCohesivePullbacks} 
Let $\calE$ be  a small extensive category with finite limits  and such that the inclusion ${R : \Dec\calE \rightarrow \calE}$ has a finite-product preserving left adjoint ${L : \calE \rightarrow \Dec\calE}$. Let ${k : \frakG(\calE) \rightarrow \frakG(\Dec\calE)}$ be the geometric morphism induced as in
Theorem~\ref{ThmGaeta}.
Let $J$ be a topology on ${\Dec\calE}$ finer than the Gaeta topology and let the following square
\[\xymatrix{
\mathcal{F} \ar[d]_-p \ar[r]  & \frakG(\calE) \ar[d]^-k \\
\Shv(\Dec\calE, J) \ar[r] & \frakG(\Dec\calE) 
}\]
be a pullback. If ${\calA(R-, \nu_C): \calA(R-,C) \rightarrow \calA(R-, R(L C))}$ is $J$-locally surjective for every $C$ in $\calE$ (where $\nu$ be the unit of ${L\dashv R}$)  then $p$ is  pre-cohesive and molecular.
\end{proposition}
\begin{proof}
The geometric morphism ${k : \frakG(\calE) \rightarrow \frakG(\Dec\calE)}$  is local, molecular and the leftmost adjoint preserves finite products by Theorem~\ref{ThmGaeta}.
Then $p$ is molecular and local  because these properties are stable under (bounded) pullback \cite[C3.3.15 and C3.6.7]{elephant}, and the leftmost adjoint preserves finite products by the same argument in Theorem~\ref{ThmGaeta}.
 Finally, as the canonical ${f_* \rightarrow f_!}$ is locally surjective  by Lemma~\ref{LemLocalSurjectivity}, the canonical ${p_* \rightarrow p_!}$ is epic by \cite[Lemma~{1.10}]{MenniExercise}.
\end{proof}

The consideration of topologies other than the atomic one on the category of connected decidable objects is motivated by recent work  with V. Marra to be discussed elsewhere.

\section{EILC toposes}
\label{SecBasicToposes}

Why are the best known examples of pre-cohesive toposes molecular?
There is a deceptively simple answer: they are essential toposes over $\Set$ and, as already observed in \cite{BarrPare80}, every essential geometric morphism over $\Set$ is molecular. In more detail, the authors of \cite{BarrPare80} say that their results become simpler over $\Sets$ because: first, all functors into $\Set$ are indexed and, second, $\Set$ is Boolean.

This cannot be the end of the story. One the one hand, it leads to the idea of sufficient conditions for molecularity based on the nature of the codomain (such as Theorem~\ref{Thm}) and, on the other, it points to the following class of toposes.

\begin{definition} A topos $\calS$ is called {\em EILC} (or {\em basic}) if every essential geometric morphism with codomain $\calS$ is molecular.
\end{definition}

The term {\em basic} is not proposed  as a serious alternative to {\em EILC}; it will be used  in this section to emphasize the idea that an EILC topos $\calS$ has some features of a topos of `discrete sets' that simplifies the study of toposes over the base $\calS$, in comparison to toposes over an arbitrary  topos. Trivially, over a basic topos, you need not worry about the distinction between essential and molecular, but there is further evidence to support the intuition that the objects of a basic topos are `discrete'. For instance, consider the following result which assumes some basic knowledge of Lawvere's dimension theory  \cite{Lawvere91}  in terms of {\em levels} (i.e. essential subtoposes) and SDG where an object $T$ is called {\em tiny} if ${(-)^T}$ has a right adjoint \cite{Yetter87}.

\begin{proposition}
If $\calS$ is a basic topos then the following hold:
\begin{enumerate}
\item (Poor dimension theory) Every level of $\calS$ is an open subtopos.
\item (Lack of infinitesimals) Every pointed tiny object in $\calS$ is terminal.
\end{enumerate}
\end{proposition}
\begin{proof}
The first item follows because every level of $\calS$ must be open by \cite[A4.5.1]{elephant}.
To prove the  second item   let ${z : 1 \rightarrow T}$ be a tiny pointed object so that  there is an essential geometric morphism ${t : \calS \rightarrow \calS}$ whose inverse image is ${(-)^T}$.
Let ${i : 1 \rightarrow T^T}$ be the transposition of the identity and let the square on the left below be a pullback. 
\[\xymatrix{
P \ar[d]_-{\pi_0} \ar[r]^-{\pi_1} & 1^T \ar[d]^-{z^T}  && P \times T \ar[d]_-{\pi_0 \times T} \ar[r]^-{\pi_1 \times T} & 1^T \times T \ar[d]^-{z^T} \ar[r]^-{ev} & 1 \ar[d]^-{z} \\
1 \ar[r]_-{i}                                  & T^T                          && 1\times T \ar[r]^-{i\times T}  \ar@(rd,ld)[rr]_-{\pi_1}                                & T^T \times T \ar[r]^-{ev} & T
}\]
Then, the square on the right above is a pullback because $t$ is molecular.
As the bottom map is an isomorphism, so is the top map.
That is, ${P \times T}$ is terminal. 
So $T$ is terminal too.
\end{proof}

It is also possible to lift Example~\ref{ExOfNonMolecularMap} to the `elementary' level by proving, for any topos $\calS$, that if ${\calS^{\rightarrow}}$ is basic then $\calS$ is degenerate.

As far as I know, the first explicit mention of EILC/basic toposes is in a May 2017 public message to the categories-list which prompted no discussion so it is not unfair to infer   that little was known on the subject. I tried to prove that basic toposes are Boolean, without success. It was only recently that the first important step to understand this class of toposes was taken in \cite{Hemelaer2022}. It turns out that a presheaf topos is basic if and only if it is Boolean \cite[Proposition~{4.4}]{Hemelaer2022}; but, to my surprise, there are many non-boolean toposes that are basic. For instance, the topos of sheaves on a Hausdorff topological space is basic, but there are more. See \cite[Theorem~{3.3}]{Hemelaer2022}.

Also, among the many possible variants in the definition of basic topos, Hemelaer identifies one that is also relevant to recall here. 
A topos $\calS$ is called {\em CILC} if,  for every geometric morphism ${f : \calE \rightarrow \calS}$ such that $f^*$ is cartesian closed, $f$  is locally connected. Every Boolean topos is CILC \cite[Theorem~{5.10}]{Hemelaer2022} an so, as explained in Remark~{5.11} loc.cit.,  Theorem~\ref{Thm} follows.

 Which are the toposes $\calS$ such that every pre-cohesive ${\calE \rightarrow \calS}$ is molecular?
 Maybe every topos has this property.

\section*{Acknowledgments}   

I would like to thank G.~Janelidze,  G.~Rosolini,  T.~Streicher and also the referee for their useful comments and suggestions.


\end{document}